\RequirePackage{fix-cm}
\documentclass[preprint]{revtex4-1}
\usepackage{graphicx}
\usepackage{amsmath,tikz,color,calc,amssymb, tikz, calc, tikz-cd, dsfont}
\usepackage{microtype}
\usepackage[bitstream-charter]{mathdesign}

\usepackage{oldgerm}
\DeclareFontFamily{U}{yswab}{}	
\DeclareFontShape{U}{yswab}{m}{n}{<-> yswab}{}

\usepackage{bm}
\usepackage[final]{showkeys}

\numberwithin{equation}{section}
\usepackage[hmargin=2.5cm,vmargin={2cm,3.5cm}]{geometry}
\allowdisplaybreaks

\newcommand{\be}{\begin{equation}}
\newcommand{\ee}{\end{equation}}

\newcommand{\beqar}{\begin{eqnarray*}}
\newcommand{\eeqar}{\end{eqnarray*}}
\newcommand{\bea}{\begin{eqnarray*}}
\newcommand{\eea}{\end{eqnarray*}}

\newcommand{\beqarl}{\begin{eqnarray}}
\newcommand{\eeqarl}{\end{eqnarray}}

\definecolor{darkgreen}{rgb}{0.0,0.5,0.0}
\definecolor{darkblue}{rgb}{0.0,0.0,0.3}
\definecolor{nicosred}{rgb}{0.65,0.1,0.1}
\definecolor{light-gray}{gray}{0.7}
\allowdisplaybreaks[1]

\RequirePackage{datetime} %comment out when arxiving or submitting
\allowdisplaybreaks
\begin{document}

\title{A Fractional Hawkes process%\thanks{Grants or other notes
%about the article that should go on the front page should be
%placed here. General acknowledgments should be placed at the end of the article.}
}
%\subtitle{Do you have a subtitle?\\ If so, write it here}

%\titlerunning{Short form of title}        % if too long for running head

\author{J. Chen}
\email{chenj60@cardiff.ac.uk}
%\noaffiliation
\affiliation{School of Mathematics, Cardiff University, UK}
\author{A. G. Hawkes}
\email{a.g.hawkes@swansea.ac.uk}
%\noaffiliation
\affiliation{School of Management, Swansea University, UK}
\author{E. Scalas}
\email{e.scalas@sussex.ac.uk}
%\noaffiliation
\affiliation{Department of Mathematics, University of Sussex, UK}
%%\authorrunning{Short form of author list} % if too long for running head
%
%\institute{J. M. Chen \at
%              School of Mathematics, Cardiff University, UK \\
%%              Tel.: +123-45-678910\\
%%              Fax: +123-45-678910\\
%              \email{chenj60@cardiff.ac.uk}           \\
%%             \emph{Present address:} of F. Author  %  if needed
%           \and
%           A. G. Hawkes \at
%             School of management, Swansea University, UK \\
%\email{a.g.hawkes@swansea.ac.uk}  \\
%\and
%E. Scalas \at 
%Department of Mathematics, University of Sussex, UK\\
%\email{e.scalas@sussex.ac.uk}}

%\date{Received: date / Accepted: date}
% The correct dates will be entered by the editor

%\maketitle

\begin{abstract}
We modify ETAS models by replacing the Pareto-like kernel proposed by Ogata with a Mittag-Leffler type kernel. Provided that the kernel decays as a power law with exponent $\beta + 1 \in (1,2]$, this replacement has the advantage that the Laplace transform of the Mittag-Leffler function is known explicitly, leading to simpler calculation of relevant quantities.
\end{abstract}

\maketitle

\section{Introduction}

\noindent In 1971, Hawkes (\cite{hawkes71,hawkes71b}) introduced a class of self-exciting processes to model contagious processes. In their simpler version, these are point processes with the following conditional intensity
$$
\lambda(t|\mathcal{H}_t) = \lim_{h \to 0} \frac{\mathbb{E}(N(t,t+h)|\mathcal{H}_t)}{h}=\lambda + \alpha \int_{-\infty}^t f(t-u) \, d N (u),
$$
where $\lambda> 0$, $N(t)$ is a Hawkes self-exciting counting process, $\mathcal{H}_t$ represents the history of the process, $\alpha$ is a branching ratio that must be smaller than 1 for stability, and $f(t)$ is a suitable kernel ($f(t)$ must be a probability density function for a positive random variable).

\noindent In 1988, Ogata (\cite{ogata88}) proposed the use of a power-law kernel for self-exciting processes of Hawkes type in order to reproduce the empirical Omori law for earthquakes. Ogata's models are also known as {\em Epidemic Type Aftershock Sequence} models or ETAS models. Within this framework, it is natural to replace Ogata's power-law kernel with a Mittag-Leffler kernel and this will be the main contribution of this chapter. We will first introduce the Mittag-Leffler distribution for positive random variables, then we will define a ``fractional'' version of Hawkes processes. Spectral properties and intensity expectation will be discussed using the fact that the Laplace transform of the one-parameter Mittag-Leffler function of argument $-t^\beta$ with $\beta \in (0,1]$ is known analytically. Finally, we will present a simple algorithm based on the thinning method by Ogata \cite{ogata81} that simulates the conditional intensity process.

\section{Mittag-Leffler distributed positive random variables}

\noindent Consider the one-parameter Mittag-Leffler function
\begin{equation}
\label{MLF}
E_\beta (z) := \sum_{n=0}^\infty \frac{z^n}{\Gamma(n \beta + 1)},
\end{equation}
with $\beta \in (0,1]$. If computed on $z=-t^\beta$ for $t \geq 0$, the Mittag-Leffler function $E_\beta (-t^\beta)$ has the meaning of survival function for a positive random variable $T$ with infinite mean. This function interpolates between a stretched exponential for small times and a power-law with index $\beta$ for large times.
Its sign-changed first derivative 
\begin{equation}
\label{pdfML}
f_\beta (t) := - \frac{d E_\beta(-t^\beta)}{d t} = t^{\beta-1} E_{\beta,\beta}(-t^\beta)
\end{equation}
is the probability density function of $T$, where $E_{\alpha,\beta}(z)$ is the two-parameter Mittag-Leffler function defined as
\begin{equation}
\label{2PMLF}
E_{\gamma,\delta}(z) := \sum_{n=0}^\infty \frac{z^n}{\Gamma(n \gamma + \delta)}. 
\end{equation}
Notice that the one-parameter Mittag-Leffler function coincides with the two parameter one with $\gamma=\beta$ and $\delta=1$.
For a suitable function $f(t)$ defined for positive $t$, let us define its Laplace transform as
$$
\tilde{f} (s) = \mathcal{L} (f(t),s) = \int_0^\infty f(t) \, \mathrm{e}^{-st} \, dt.
$$
The functions $E_\beta(-t^\beta)$ and $f_\beta (t)$ have explicit Laplace transforms. The survival function has the following Laplace transform
\begin{equation}
\label{LTS1}
\mathcal{L}(E_\beta (-t^\beta); s) = \frac{s^{\beta -1}}{1+s^\beta},
\end{equation}
and the probability density function has the following Laplace transform
\begin{equation}
\label{LTS2}
\mathcal{L}(f_\beta(t); s) = \frac{1}{1+s^\beta}.
\end{equation}
Moreover, they have an explicit representation as an infinite (actually continuous) sum of exponential functions \cite{mainardigorenfloln};
\begin{equation}
\label{sum1}
E_\beta (-t^\beta) = \int_0^\infty \mathrm{e}^{-\theta t} K_\beta (\theta) \, d \theta,
\end{equation}
with
\begin{equation}
\label{sumkernel}
K_\beta (\theta) = \frac{1}{\pi} \frac{\theta^{\beta-1} \sin(\beta \pi)}{\theta^{2\beta}+2 \theta^\beta \cos(\beta \pi)+1},
\end{equation}
leading to
\begin{equation}
\label{sum2}
f_\beta (t) = \int_0^\infty \theta \mathrm{e}^{-\theta t} K_\beta (\theta) \, d \theta.
\end{equation}

\noindent These functions play an important role in fractional calculus. For instance $E_\beta (-t^\beta)$ is the solution of the following anomalous relaxation problem
\begin{equation}
\label{relax}
\frac{d^\beta g(t)}{d t^\beta} = - g(t),
\end{equation}
where $d^\beta / d t^\beta$ is the Caputo derivative defined as
\begin{equation}
\label{Caputo}
\frac{d^\beta g(t)}{d t^\beta} = \frac{1}{\Gamma(1-\beta)} \frac{d}{dt} \int_0^t \frac{g(\tau)}{(t-\tau)^\beta} \, dt -
\frac{t^{-\beta}}{\Gamma(1-\beta)} g(0^+),
\end{equation}
with initial condition $g(0^+)=1$.

\section{The fractional Hawkes processes}

It becomes natural to use $f_\beta (t)$ as kernel for a version of Hawkes processes that we can call {\em fractional} Hawkes processes. 
\noindent The conditional intensity of the process is given by
\begin{equation}
\label{fHp}
\lambda(t|\mathcal{H}_t) = \lim_{h \to 0} \frac{\mathbb{E}(N(t,t+h)|\mathcal{H}_t)}{h}=\lambda + \alpha \int_{-\infty}^t f_\beta(t-u) \, d N (u),
\end{equation}
where $\lambda> 0$ and $N(t)$ is a Hawkes self-exciting point process, leading to
\begin{equation}
\label{fHp1}
\lambda(t|\mathcal{H}_t) =\lambda + \alpha \sum_{t_i < t}  f_\beta (t-t_i).
\end{equation}
with the branching ratio $\alpha<1$ for stability. Hainaut \cite{hainault19} gives a different definition of fractional Hawkes process. Let us use his notation in the following remark. In his paper, he considers the time-changed intensity process $\lambda_{S_t}$ where, in his case, the conditional intensity $\lambda_t$ of the self-exciting process is the solution of a mean-reverting stochastic differential equation
$$
d \lambda_t = \kappa (\theta - \lambda_t) \, dt + \eta dP_t,
$$
where $\kappa$, $\theta$ and $\eta$ are suitable parameters and the driving process $P_t$ is given by
$$
P_t := \sum_{k=1}^{N_t} \xi_k,
$$
where $N_t$ is a counting process and $\xi_i$s are independent and identically distributed marks with finite positive mean and finite variance. The time-change $S_t$ is the inverse of a $\beta$-stable subordinator. Our definition \eqref{fHp} is much simpler and it is directly connected with ETAS processes given that the kernel $f_\beta (t)$ has power-law tail with index $\beta+1$. \cite{mainardi00}. In particular, given the explicit Laplace transform of $f_\beta (t)$ and its representation in terms of infinite sum of exponentials, we can derive some explicit formulas.

\subsection{Spectral properties}

From equation (11) in \cite{hawkes71}, we get the following equation for the covariance density $\mu(\tau)$ for $\tau >0$.
\begin{equation}
\label{covdens}
\mu(\tau) = \alpha \left[\Lambda f_\beta(\tau) + \int_0^\tau f_\beta (\tau - v) \, \mu(v) \, dv + \int_0^\infty f_\beta(\tau+v) \, \mu(v) \, dv \right],
\end{equation}
where $\Lambda$ represents the asymptotic stationary value of the conditional intensity as derived in equation (3.9) below.
If we now take the Laplace transform of \eqref{covdens}, we get
\begin{equation}
\label{covdenslp}
\tilde{\mu}(s) = \alpha \left[\Lambda \tilde{f}_\beta(s) + \tilde{f}_\beta(s)  \tilde{\mu} (s) + \mathcal{L} \left( \int_0^\infty f_\beta(\tau+v) \, \mu(v) \, dv; s \right) \right].
\end{equation}
Now, using \eqref{sum2} and setting
$$
h(\theta) = \theta K_\beta(\theta) =  \frac{1}{\pi} \frac{\theta^{\beta} \sin(\beta \pi)}{\theta^{2\beta}+2 \theta^\beta \cos(\beta \pi)+1},
$$
we can write
$$
f_\beta (\tau) = \int_0^\infty h(\theta) \mathrm{e}^{-\theta \tau} \, d \theta,
$$
so that the last summand in \eqref{covdenslp} becomes
\begin{eqnarray}
\label{lastsum}
\mathcal{L} \left( \int_0^\infty f_\beta(\tau+v) \, \mu(v) \, dv; s \right) & = & \int_0^\infty \mathrm{e}^{-s \tau}
\left[\int_0^\infty \left[\int_0^\infty h(\theta) \mathrm{e}^{-\theta(\tau+v)} \, d\theta \right] \mu(v) \, dv \right] \, d\tau \nonumber \\
& = & \int_0^\infty h(\theta) \frac{1}{\theta+s} \tilde{\mu} (\theta) \, d \theta.
\end{eqnarray}
In principle, a numerical approximation of the integral in equation \eqref{lastsum} plugged into \eqref{covdenslp} and coupled with equation \eqref{LTS2} can lead to an explicit approximate expression for the Laplace transform $\tilde{\mu} (s)$. This will be the subject of a further paper. An alternative approach is given in \cite{hawkes71b} where the Bartlett spectrum is defined for real $\omega$ as
\begin{equation}
\label{bsp}
f(\omega) = \frac{1}{2 \pi} \int \mathrm{e}^{-i \omega \tau} \mu^{(c)} (\tau) \, d \tau,
\end{equation}
where, because $\mathbb{E} [(dN(t))^2] = \mathbb{E}[d N(t)]$ if events cannot occur multiply, the complete covariance density contains a delta function
$$
\mu^{(c)} (\tau) = \Lambda \delta (t) + \mu (t).
$$
Then it is shown in \cite{hawkes71b} p. 441 that
\begin{equation}
\label{bspf}
f(\omega) = \frac{\Lambda}{2 \pi (1-G(\omega))(1-G(-\omega))},
\end{equation}
where, in our case, we would have
$$
G(\omega) = \int_0^\infty \mathrm{e}^{-i \omega \tau} \alpha f_\beta (\tau) \, d \tau = \frac{\alpha}{1+(i \omega)^\beta}.
$$
The proof of this result depended on the assumption that the exciting kernel decays exponentially asymptotically. However, this is not true for the Mittag-Leffler distribution, which decays as a power law. Bacry and Muzy \cite{bacry}
prove a more general result using Laplace transforms in the complex plane, more easily digested from section 2.3.1 in \cite{bacrybis}. Then the Laplace transform of the covariance density is given by
\begin{equation}
\label{bslt}
\tilde{\mu}^{(c)} (s) = \frac{\Lambda}{(1 - \tilde{\Phi}(s))(1-\tilde{\Phi}(-s))},
\end{equation}
where
$$
\tilde{\Phi}(s) = \int_0^\infty \mathrm{e}^{-s \tau} \Phi(\tau) \, d \tau
$$
is the Laplace transform of the exciting kernel. In our case for $\tau >0$, we have
$$
\Phi(\tau) = \alpha f_\beta (\tau)
$$
and
$$
\tilde{\Phi}(s) = \frac{\alpha}{1+s^\beta}.
$$
Equations \eqref{bspf} and \eqref{bslt} look much the same, apart from a change of notation and a factor $2 \pi$. The difference, however, is that in \eqref{bspf} we are dealing with real $\omega$ while, in \eqref{bslt}, $s$ is a general complex variable and we can choose its domain to obtain well-behaved functions.

\subsection{Intensity expectation}

Let us consider the expectation $\Lambda (t) = \mathbb{E}[\lambda(t|\mathcal{H}_t)]$ for both a stationary and non-stationary fractional Hawkes process.
In the stationary case (process from $t=-\infty$), from equation \eqref{fHp}, we get $\Lambda = \lambda + n \Lambda$, leading to
\begin{equation}
\label{stationarymean}
\Lambda = \frac{\lambda}{1-\alpha}.
\end{equation}
On the contrary in the non-stationary case (process from $t=0$), we can modify equation \eqref{fHp} as follows
\begin{equation}
\label{fHpns}
\lambda(t|\mathcal{H}_t) =\lambda + \alpha \int_{0}^t f_\beta(t-u) \, d N (u),
\end{equation}
so that the time-dependent expectation obeys the equation
\begin{equation}
\label{fHpns1}
\Lambda(t) =\lambda + \alpha \int_{0}^t f_\beta(t-u) \, \Lambda(u) \,du.
\end{equation}
Taking Laplace transforms, we get
$$
\tilde{\Lambda} (s) = \frac{\lambda}{s} + \alpha \tilde{f}_\beta (s) \tilde{\Lambda} (s),
$$
so that
$$
\tilde{\Lambda} (s) = \frac{\lambda}{s} \frac{1}{1 -\alpha \tilde{f}_\beta (s)}.
$$
Using equation \eqref{LTS2} yields
\begin{equation}
\label{nonstationarymean}
\tilde{\Lambda} (s) = \frac{\lambda}{s} \frac{1+s^\beta}{(1-\alpha)+s^\beta}.
\end{equation}
Equation \eqref{nonstationarymean} can be inverted numerically (or analytically for $\beta =1/2$) to give $\Lambda (t)$ as well as the expected number of events from $0$ to $t$ as 
$$
\mathbb{E}[N(t)]= \int_0^t \Lambda(\tau) \, d \tau.
$$
Also, based on a continuous version of Hardy-Littlewood Tauberian theorem \cite{feller}
we get
$$
\lim_{t \to \infty} \Lambda(t) = \frac{\lambda}{1-\alpha}
$$
as given by equation \eqref{stationarymean}. This result is exemplified in Fig. 1 for $\beta = 1/2$, $\lambda = 1$ and $\alpha=1/2$. In that case, we get, for $t>0$
$$
\Lambda(t) = 2 - \mathrm{e}^{t/4} \mathrm{erfc}( \sqrt{t}/2 ).
$$
\begin{figure}
\label{tauberian}
\includegraphics[width=0.9\textwidth]{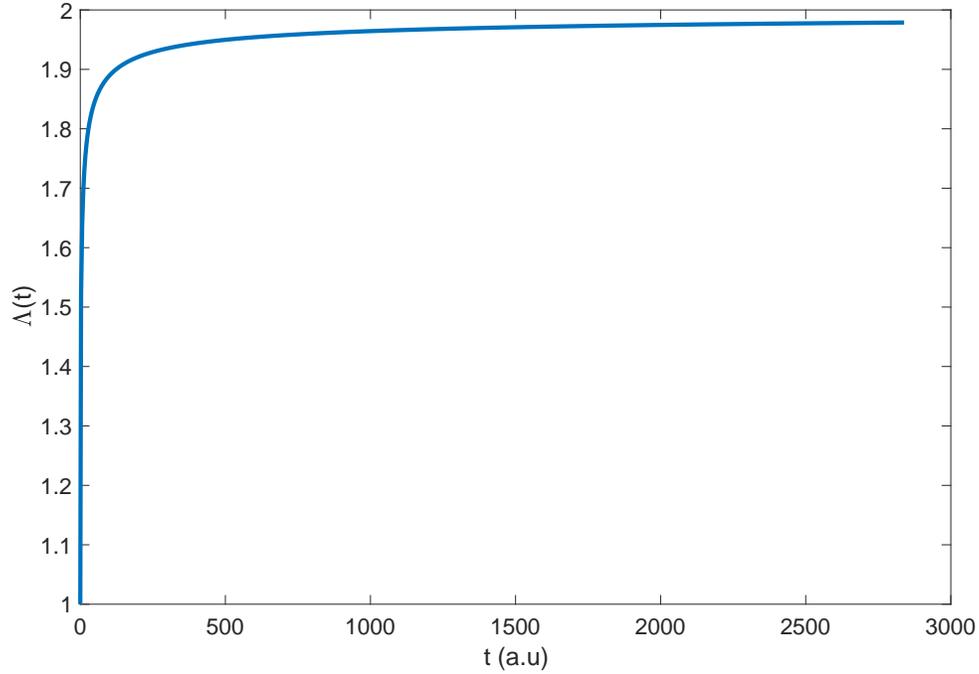}
\caption{$\Lambda(t)$ as a function of $t$ for $\beta = 1/2$, $\lambda = 1$ and $\alpha=1/2$.}
\end{figure}
From Fig. 1, one can see that $\Lambda(t)$ goes up very fast at first and, then, slowly converges to its asymptotic value $2$. This is presumably due to the long tail of the Mittag-Leffler kernel.

\section{Simulation}

In order to simulate the intensity process introduced in equation \eqref{fHp}, we use the thinning algorithm introduced by Ogata \cite{ogata81} (see also Zhuang and Touati 2015 report \cite{zhuang15}). The function {\tt ml.m} described in \cite{garrappa} is needed to compute the Mittag-Leffler functions described above and can be retrieved from the Matlab file exchange. The algorithm is as follows
\begin{enumerate}
\item Set the initial time $t = 0$, a counter $i = 0$ and a final time $T$.
\item Compute $ M= \lambda + \alpha \sum_{t_i < t+\varepsilon}  f_\beta (t+\varepsilon-t_i)$, for some small value of $\varepsilon$.
\item Generate a positive exponentially distributed random variable $E$ with the meaning of a waiting time, with rate $1/M$.
\item Set $\tau=t+E$.
\item Generate a uniform random variate $U$ between $0$ and $1$.
\item If $U < [\lambda + \alpha \sum_{t_i < \tau}  f_\beta (\tau-t_i)]/M$, set $t_{i+1} = \tau$  and update time to $t=\tau$, else, just set $t=\tau$.
\item Return to step 2 until $t$ exceeds $T$.
\item Return the set of event times (or epochs) $t_i$
\end{enumerate}
An implementation of this algorithm in Matlab is presented below. \\

{\tt 
\noindent T=5; \\
\noindent t=0; \\
\noindent n=0; \\
alpha=0.9; \% For stability \\
\noindent mu=1; \\
\noindent epsilon=1e-10; \\
\noindent beta=0.7; \\

\noindent SimPoints=[]; \\

\noindent while t<T \\
\noindent     t\\
\noindent    M=mu+sum(alpha*ml(-(t+epsilon-SimPoints).$^\wedge$beta,beta,beta,1). \dots\\
\noindent*(t+epsilon-SimPoints).$^\wedge$(beta-1)); \\
\noindent    E=exprnd(1/M,1,1); \\
\noindent    t=t+E; \\
\noindent   U=rand;\\
\noindent    if (U<(((mu+sum(alpha*ml(-(t-SimPoints).$^\wedge$beta,beta,beta,1). \dots\\
\noindent*(t-SimPoints).$^\wedge$(beta-1))))/M)) \\
\noindent        n=n+1;\\
\noindent        SimPoints = [SimPoints, t];\\
\noindent    end\\
\noindent end\\
\noindent index=find(SimPoints<10);\\
\noindent SimPoints=SimPoints(index);\\
}

\noindent Two examples of intensity process simulated up to $t=5$ are presented in Fig. 2 and Fig. 3 for $\beta = 0.9$ and $\beta=0.7$, respectively. 
\begin{figure}
\label{intensity9}
\includegraphics[width=0.9\textwidth]{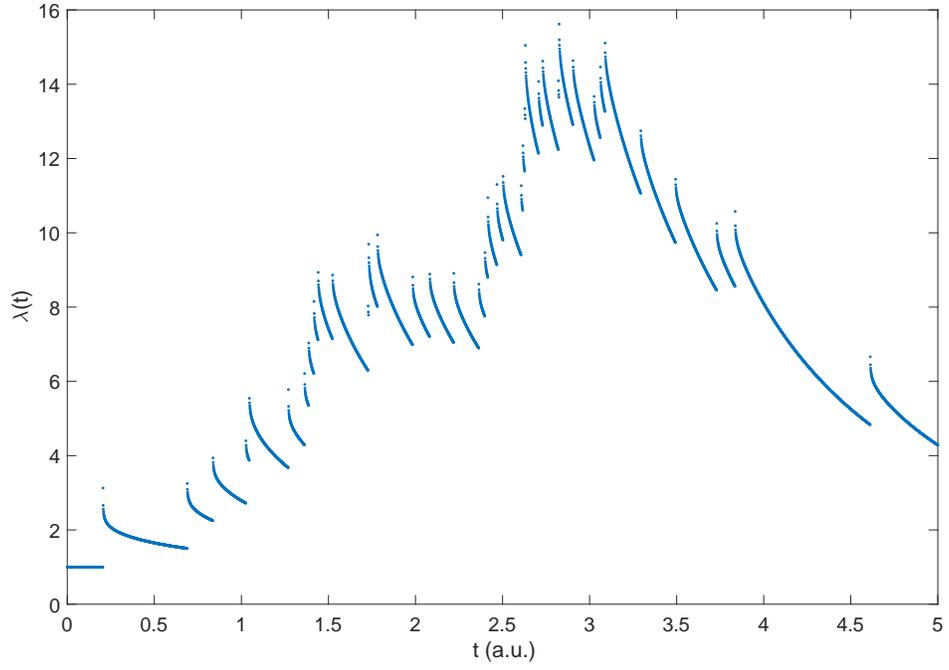}
\caption{$\lambda(t|\mathcal{H}_t)$ as a function of time $t$ for $\beta = 0.9$, $\lambda = 1$ and $\alpha=0.9$.}
\end{figure}
\begin{figure}
\label{tauberian}
\includegraphics[width=0.9\textwidth]{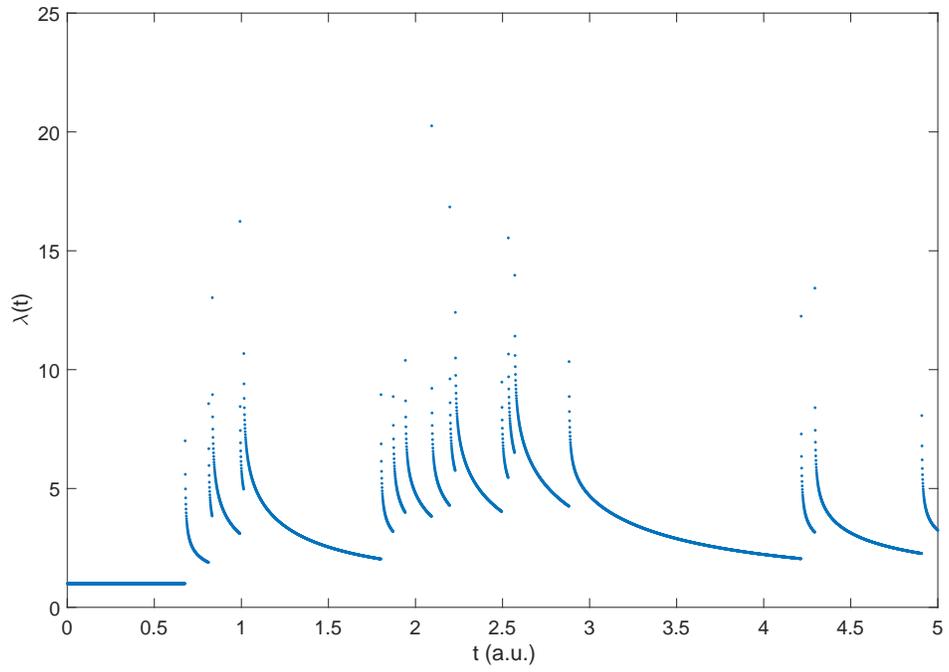}
\caption{$\lambda(t|\mathcal{H}_t)$ as a function of $t$ for $\beta = 0.7$, $\lambda = 1$ and $\alpha=0.9$.}
\end{figure}

\section{Outlook}

In this paper, we defined a ``fractional'' Hawkes process and we studied its spectral properties and the expectation of its intensity. Thanks to explicit expressions for Laplace transforms, we could derive some analytical expressions that are only asymptotically available for power-law kernels of Pareto type as originally suggested by Ogata. We also presented an explicit simulation of the intensity process based on the so-called thinning method.

Further work is needed to better characterize our process. In particular, we did not deal with parameter estimation, the multivariate version of the process, and we did not use the model to fit earthquake data (or any other data for what matters including financial data). We do hope that all this and more can become the subject of an extensive future paper on this process.

%\begin{acknowledgements}
%If you'd like to thank anyone, place your comments here
%and remove the percent signs.
%\end{acknowledgements}

% BibTeX users please use one of
%\bibliographystyle{spbasic}      % basic style, author-year citations
%\bibliographystyle{spmpsci}      % mathematics and physical sciences
%\bibliographystyle{spphys}       % APS-like style for physics
%\bibliography{}   % name your BibTeX data base

% Non-BibTeX users please use

\end{document}